\documentclass[12pt]{amsart}

\usepackage{amssymb}
\usepackage{amsfonts}
\usepackage{amsmath}

\theoremstyle{plain}
\newtheorem*{thm*}{Theorem}
\newtheorem*{prop*}{Proposition}

\begin{document}
\baselineskip=17pt


\title[Continuity of Seminorms on Finite-Dimensional Vector Spaces]{Continuity of Seminorms on \\ Finite-Dimensional Vector Spaces}

\author{Moshe Goldberg}

\address{Department of Mathematics,
Technion -- Israel Institute of Technology,
Haifa 32000, Israel}

\email{mg@technion.ac.il}



\begin{abstract}
The main purpose of this note is to establish the continuity of seminorms on finite-dimensional
vector spaces over the real or complex numbers.
\end{abstract}

\maketitle

Throughout this note let $\mathbf{V}$ be a finite-dimensional vector space over a field $\mathbb{F}$, either $\mathbb{R}$ or $\mathbb{C}$. Let $N$ be a \textit{norm} on $\mathbf{V}$, so that for all
$a,b \in \mathbf{V}$ and $\alpha \in \mathbb{F}$,
\begin{equation*}
\begin{split}
& N(a)>0, \quad a \neq 0, \\
& N(\alpha a) = |\alpha|N(a), \\
& N(a+b) \leq N(a) + N(b).
\end{split}
\end{equation*}
Since $\mathbf{V}$ is finite-dimensional, all norms on $\mathbf{V}$ are equivalent, inducing on $\mathbf{V}$ a unique topology. In particular, it follows that all norms on $\mathbf{V}$ are continuous with respect to this topology.

The notion of norm can be relaxed in two familiar ways:
\vspace{6pt}

(a) We say that a real-valued function $f: \mathbf{V} \rightarrow \mathbb{R}$ is a \textit{subnorm} on $\mathbf{V}$ if for all $a \in \mathbf{V}$ and $\alpha \in \mathbb{F}$,
\begin{equation*}
\begin{split}
& f(a) > 0, \quad a \neq 0, \\
& f(\alpha a) = |\alpha|f(a).
\end{split}
\end{equation*}

(b) We say that a real-valued function $S: \mathbf{V} \rightarrow \mathbb{R}$ is a \textit{seminorm} on $\mathbf{V}$ if for all $a,b \in \mathbf{V}$ and $\alpha \in \mathbb{F}$,
\begin{equation*}
\begin{split}
& S(a) \geq 0, \\
& S(\alpha a) = |\alpha|S(a), \\
& S(a+b) \leq S(a) + S(b).
\end{split}
\end{equation*}

It follows that a subnorm $f$ is a norm on $\mathbf{V}$ if and only if $f$ is \textit{subadditive}. Similarly, a seminorm $S$ is a norm on $\mathbf{V}$ if and only if $S$ is \textit{positive-definite}.

If $\dim \mathbf{V} = 1$, then fixing a nonzero element $a_{0} \in \mathbf{V}$, we may write
$$
\mathbf{V} = \{ \alpha a_{0}: \alpha \in \mathbb{F} \}.
$$
So every subnorm $f$ on $\mathbf{V}$ must be of the form
$$
f(a) = \gamma |\alpha|, \quad a = \alpha a_{0} \in \mathbf{V},
$$
where $\gamma$, the value of $f$ at $a_{0}$, is a positive constant and, consequently, $f$ is continuous
on $\mathbf{V}$.

If, however, $\dim \mathbf{V} \geq 2$, then contrary to norms, subnorms on $\mathbf{V}$ may fail to be continuous. For example, \cite[Section 3]{G}, let $f$ be a continuous subnorm on $\mathbf{V}$. Fix an element $a_{0} \neq 0$ in $\mathbf{V}$, and let
$$
\mathbf{W} = \{ \alpha a_{0} : \alpha \in \mathbb{F} \}
$$
be the one-dimensional linear subspace of $\mathbf{V}$ generated by $a_{0}$. Select a real number $\kappa$, $\kappa > 1$, and set
\begin{equation*}
g_{\kappa} =
\begin{cases}
\kappa f(x),  & a \in \mathbf{W}, \\
f(x), & a \in \mathbf{V} \smallsetminus \mathbf{W}.
\end{cases}
\end{equation*}
Then, evidently, $g_\kappa$ is a subnorm on $\mathcal{A}$. Further,
$g_\kappa$ is discontinuous at $a_0$ since
$$
\lim_{\substack{a \to a_0 \\ a \notin \mathbf{W}}} g_\kappa(a)
= \lim_{\substack{a \to a_0 \\ a \notin \mathbf{W}}} f(a)
= \lim_{a \to a_0} f(a)=f(a_0) \neq g_\kappa(a_0).
$$

In fact, subnorms on finite-dimensional vector spaces can be discontinuous \textit{everywhere}. Examples of such pathological subnorms on the complex numbers and on the quaternions (taken as 2- and 4-dimensional vector spaces over the reals), have emerged in parts (a) and (c) of Theorems 2.1 and 3.1 of \cite{GL2}.

As a last introductory remark, we note that just like norms, all continuous subnorms on a finite-dimensional vector space $\mathbf{V}$ are equivalent to each other (e.g., \cite[Lemma 1.1]{GL1}); i.e., if $f$ and $g$ are continuous subnorms on $\mathbf{V}$, then there exist positive constants $\mu > 0$,
$\nu > 0$, such that
$$
\mu f(a) \leq g(a) \leq \nu f(a), \quad a\in \mathbf{V}.
$$

With these preliminaries, we turn now to discuss the continuity of seminorms.

\begin{thm*}
Let $S$ be a seminorm on a finite-dimensional vector space $\mathbf{V}$ over $\mathbb{F}$, either $\mathbb{R}$ or $\mathbb{C}$. Then $S$ is continuous with respect to the unique topology on $\mathbf{V}$.
\end{thm*}

We shall provide two different proofs:

\begin{proof}[Proof 1]
Let $\{ e_{1}, \ldots, e_{n}\}$ be a basis for $\mathbf{\mathbf{V}}$, so that every $a \in \mathbf{V}$ can be uniquely expressed as
$$
a = \alpha_{1}e_{1} + \cdots + \alpha_{n}e_{n}, \quad (\alpha_{j} \in\ \mathbb{F}).
$$
If $S=0$, then there is nothing to prove, so we may assume that $S \neq 0$ which implies that
$$
\sigma = \max_{1 \leq j \leq n} S(e_{j})
$$
is a positive constant. Hence,
\begin{equation}
\label{eq:1}
N_{1}(a) = \sigma \sum_{j=1}^{n} |\alpha_{j}|, \quad
a = \alpha_{1}e_{1} + \cdots + \alpha_{n}e_{n} \in \mathbf{V},
\end{equation}
is a norm on $\mathbf{V}$; moreover, $N_{1}$ majorizes $S$ on $\mathbf{V}$ since
\begin{equation}
\label{eq:2}
S(a) = S(\alpha_{1}e_{1} + \cdots + \alpha_{n}e_{n}) \leq
|\alpha_{1}|S(e_{1}) + \cdots + |\alpha_{n}|S(e_{n}) \leq
\sigma \sum_{j=1}^{n} |\alpha_{j}| = N_{1}(a).
\end{equation}

Further, for all $a,b \in \mathbf{V}$, we get
\begin{subequations}
\label{3}
\begin{gather}
\label{3a}
S(a) - S(b) = S(a-b+b) - S(b) \leq S(a-b) + S(b) - S(b) = S(a-b), \\
\intertext{as well as}
\label{3b}
S(b) - S(a) = S(b-a+a) - S(a) \leq S(b-a) + S(a) - S(a) = S(a-b).
\end{gather}
\end{subequations}
Thus, by (2) and (3),
$$
|S(a) - S(b)| \leq N_{1}(a-b).
$$
It follows that if $\{ a_{j} \}_{j=1}^{\infty}$ is a sequence in $\mathbf{V}$ which tends to a vector $b$, then
$$
|S(a_{j})-S(b)| \leq N_{1}(a_{j}-b)\rightarrow 0 \quad \textrm{as } j \rightarrow \infty,
$$
and we are done.
\end{proof}

\begin{proof}[Proof 2]
Since $S$ is a seminorm, $\mathcal{K} = \ker S$ is a subspace of $\mathbf{V}$. Consider the quotient space $\mathbf{V}/\mathcal{K}$, and set
\begin{equation}
\label{eq:4}
N_{2}(a+\mathcal{K}) = S(a), \quad a\in\mathbf{V}.
\end{equation}
We note that if $a+\mathcal{K} = b+\mathcal{K}$ for some $a$ and $b$ in $\mathbf{V}$, then
$a-b \in \mathcal{K}$; so
$$
S(a) = S(a-b+b) \leq S(a-b) + S(b) = S(b),
$$
and simiarly, $S(b) \leq S(a)$. Hence $S(a) = S(b)$, and it follows that $N_{2}$ is well defined on
$\mathbf{V}/\mathcal{K}$. Furthermore, it is easy to verify that $N_{2}$ is a norm on $\mathbf{V}/\mathcal{K}$.

Now, since $\mathbf{V}$ is finite dimensional, so is $\mathbf{V}/\mathcal{K}$ and, consequently, $N_{2}$ is continuous on $\mathbf{V}/\mathcal{K}$. In addition, if $\{ a_{j} \}_{j=1}^{\infty} \subset \mathbf{V}$ is a sequence that tends to $b$, then $a_{j} + \mathcal{K} \rightarrow b + \mathcal{K}$ since the mapping
$a \rightarrow a + \mathcal{K}$ is linear and therefore continuous. Whence,
$$
S(a_{j})-S(b) = N_{2}(a_{j}+\mathcal{K}) - N_{2}(b+\mathcal{K}) \rightarrow 0 \quad
\textrm{as } j \rightarrow \infty,
$$
and the proof is complete.
\end{proof}

In passing, we observe that for every $a \in \mathbf{V}$ and $b \in \mathcal{K}$,
$$
S(a) = S(a-b+b) \leq S(a-b) + S(b) = S(a-b) \leq S(a) + S(b) = S(a).
$$
Hence, by (4),
$$
N_{2}(a+\mathcal{K}) = S(a-b), \quad a\in\mathbf{V},~b\in\mathcal{K};
$$
so $N_{2}(a+\mathcal{K})$ \textit{is the distance, as measured by} $S$, \textit{from any} $a \in \mathbf{V}$ \textit{to any} $b \in \mathcal{K} $.

We conclude this short note by posting another simple result.

\begin{prop*}
Let $\mathbf{V}$ be a finite-dimensional vector space over $\mathbb{F}$, equipped with a seminorm $S$ and a norm $N$. Then $S$ is left-equivalent to $N$, i.e., there exists a constant $\tau>0$ such that
$$
S(a) \leq \tau N(a) \quad \textit{for all } a\in\mathbf{V}.
$$
\end{prop*}

As for our theorem, we provide two proofs.

\begin{proof}[Proof 1]
Since $N_{1}$ in (1) is equivalent to $N$, we can find a positive constant $\tau>0$ for which
$$
N_{1}(a) \leq \tau N(a), \quad a\in\mathbf{V}.
$$
Hence, by (2), the desired result follows.
\end{proof}

\begin{proof}[Proof 2]
We may assume that $S$ is not identically zero. Consider the compact set
$$
\mathcal{B} = \{a\in\mathbf{V}: ~ N(a)=1\}.
$$
By our theorem, $S$ is continuous on $\mathcal{B}$; and since $S$ may not vanish identically on $\mathcal{B}$, we infer that
$$
\tau = \max\{S(a): ~ a\in\mathcal{B}\}
$$
is a positive constant. Now, select any $a\in\mathbf{V}$, $a\neq0$. Then $a/N(a)\in\mathcal{B}$; so
$$
S\left( \frac{a}{N(a)} \right) \leq \tau,
$$
and the assertion is in the bag.
\end{proof}

\end{document}